\documentclass[11pt]{article}
\usepackage[a4paper,margin=1in]{geometry}

\usepackage[T1]{fontenc}

\usepackage{hyperref}
\usepackage{indentfirst}
\usepackage{booktabs}
\usepackage{color}

\usepackage{amsmath}
\usepackage{amssymb}

\usepackage{float} 

\usepackage[utf8]{inputenc} 

\usepackage[normalem]{ulem}

\usepackage{booktabs}
\usepackage{float}
\usepackage{graphicx}
\usepackage{caption}

\usepackage{amsthm}

\newtheorem{theorem}{Theorem}
\newtheorem{definition}{Definition}

\usepackage{filecontents}
\usepackage{xr}
\usepackage{hyperref}

\DeclareUnicodeCharacter{FF0C}{,}

\newcommand{\Xt}{\{X_t\}_{t \ge 0}}

\usepackage{fancyhdr}
\usepackage{lipsum}  

\newcommand{\authoraffil}{\noindent
$^{1}$ Department of Mathematics, Tokyo University of Science, Chiba, Japan , \texttt{6124507@ed.tus.ac.jp}\\
$^{2}$ Department of Industrial and Systems Engineering, \texttt{ishigaki@rs.tus.ac.jp}
}

\fancypagestyle{firstpagefooter}{
  \fancyhf{} 
   
  \fancyfoot[L]{\scriptsize \authoraffil}
  \fancyfoot[R]{\thepage}
}

\title{Inventory Control Using a Lévy Process for Evaluating Total Costs under Intermittent Demand}

\author{
Ryoya Koide$^{1}$, Yurika Ono$^{2}$, Aya Ishigaki$^{2}$ \\
}
\date{ \ }

\begin{document}
\maketitle              
\thispagestyle{firstpagefooter} 
\begin{abstract}
Products with intermittent demand are characterized by a high risk of sales losses and obsolescence due to the sporadic occurrence of demand events. Generally, both point forecasting and probabilistic forecasting approaches are applied to intermittent demand. In particular, probabilistic forecasting, which models demand as a stochastic process, is capable of capturing uncertainty. An example of such modeling is the use of Lévy processes, which possess independent increments and accommodate discontinuous changes (jumps). However, to the best of our knowledge, in inventory control using Lévy processes, no studies have investigated how the order quantity and reorder point affect the total cost. One major difficulty has been the mathematical formulation of inventory replenishment triggered at reorder points.
To address this challenge, the present study formulates a reorder-point policy by modeling cumulative demand as a drifted Poisson process and introducing a stopping time to represent the timing at which the reorder point is reached. Furthermore, the validity of the proposed method is verified by comparing the total cost with that obtained from a case where an ARIMA model is combined with a reorder-point policy.
As a main result, while the total cost under ARIMA-based forecasting increases linearly over time, the Lévy process-based formulation provides an analytical expression for the total cost, revealing that random demand fluctuations cause the expected total cost to grow at a rate faster than linear.
     
\end{abstract}
\section{Introduction}
\noindent In manufacturing settings, supply and demand operations play a critical role in preventing overstocking and stockouts of products or raw materials, thereby reducing costs and enhancing customer satisfaction. 
Decision-makers in supply and demand operations determine order quantities and adjust inventory levels based on demand forecasts. However, certain types of demand data are difficult to predict. 
A notable example is intermittent demand, which consists of count data characterized by a significant number of zero observations~\cite{WANG20241038}. 
Intermittent demand is commonly observed in the spare parts used during the production phases of industries such as automobiles, aircraft, and electronic equipment, accounting for approximately 60\% of inventory~\cite{KOURENTZES2013198}. 
Furthermore, spare parts are subject to after-sales demand and must respond to customer needs even after the product purchase, leading to highly irregular demand fluctuations and increasing the risks of stockouts and obsolescence~\cite{yuna2023inventory}. In many cases, handling intermittent demand relies heavily on the intuition and experience of skilled workers, creating a need for measures to reduce dependency on individual expertise in inventory management.
\par Intermittent demand is characterized by the uncertainty of both nonzero demand sizes and inter-arrival intervals. In general, demand forecasting for intermittent demand includes point forecasting and probabilistic forecasting. However, few studies have focused on probabilistic forecasting combined with the evaluation of inventory performance, leaving this area relatively unexplored~\cite{WANG20241038}.
In existing studies adopting probabilistic forecasting, intermittent demand has been modeled using stochastic processes such as Poisson or compound Poisson processes. For instance, Noba and Yamasaki~\cite{MR4702096} considered inventory control based on expected utility by employing Lévy processes. A Lévy process is a type of stochastic process that accommodates discontinuous changes (jumps) while maintaining the property of independent increments, offering strong theoretical advantages. Here, expected utility refers to the expected level of satisfaction derived from decision-making under uncertainty. In their framework, expected utility is calculated using scale functions, which are Laplace-transformed functions that enable tractable numerical analysis even for general Lévy processes.
\par However, to the best of our knowledge, no study has concretely examined how each cost component (holding cost, shortage cost, and ordering cost), as well as the reorder point and order quantity, affects the total cost when intermittent demand is modeled as a stochastic process.
\par To address this challenge, this study makes the following contributions: 
(1) The cumulative demand is formulated as a drifted Poisson process, defined as the sum of a linear function and a Poisson process. The linear component represents regular daily demand, whereas the Poisson process captures abrupt increases in demand.
(2) The concept of a \emph{stopping time} is introduced to appropriately model the timing of reorder points, thereby enabling the formulation of inventory replenishment policies.
(3) The validity of the proposed approach is verified by comparing the total cost obtained through the stochastic process-based inventory control with that resulting from the combination of ARIMA-based demand forecasting and a standard reorder point policy.
\par The contributions of this study are summarized as follows:
(1) By formulating cumulative demand using a specific Lévy process (the drifted Poisson process) and representing the expected inventory level at any given time in an intuitively interpretable form, this study examines the effects of major parameters—such as ordering cost, shortage cost, holding cost, order quantity, and ordering policies—on the total inventory cost under intermittent demand.
(2) Inventory is defined as a random variable given by "initial inventory" $-$ "cumulative demand" $+$ "cumulative replenishment," which enables elementary-level calculation of expected inventory levels.
(3) By formulating uncertainty through a stochastic process, the information contained in the demand data is preserved during the analysis, allowing for an investigation into how randomness influences the expected total cost.
The remainder of this paper is organized as follows: the next section introduces the stochastic model, followed by the formulation of the expected total cost. This studythen present the simulation results, compare them with those obtained using the ARIMA model with a reorder point policy, and conclude the paper with a discussion.

\section{Literature review}
\subsection{Prior Research on Forecasting and Inventory Management for Intermittent Demand}

\noindent Forecasting methods for intermittent demand can be broadly classified into two categories: point forecasting and probabilistic forecasting~\cite{WANG20241038}. 
Examples of point forecasting methods include Croston’s method~\cite{croston1972forecasting}, which separately forecasts nonzero demand sizes and inter-arrival intervals, and studies that adapt machine learning techniques such as neural networks and Holt’s method for intermittent demand forecasting~\cite{KOURENTZES2013198},~\cite{altay2008adapting}.
\par However, these point forecasting methods are limited in that they are based on a single value, such as the mean or median, and therefore cannot fully capture uncertainty, which may be inadequate for decision-making in inventory management~\cite{WANG20241038}. In contrast, probabilistic forecasting approaches model intermittent demand using stochastic processes such as Poisson or compound Poisson processes, allowing for the consideration of uncertainty. 
Nevertheless, Wang et al.~\cite{WANG20241038}, who proposed a combined probabilistic forecasting approach for intermittent demand, showed that 
the forecasting method with the highest accuracy does not necessarily result in the lowest inventory cost.
Thus, it is crucial to evaluate overall performance, including both demand forecasting and inventory management.
\par In recent years, deep learning-based forecasting methods such as recurrent neural networks (RNNs), 
long short-term memory networks (LSTM), and Trans \-former-based architectures have also been applied to intermittent demand prediction~\cite{Zhang2003},\cite{KOURENTZES2013198}, \cite{Zhou2021}. 
These methods often achieve high predictive accuracy. However, they tend to operate as black-box models and lack analytical interpretability, making them less suitable for downstream tasks such as cost-based inventory control.
\par To address these limitations, the present study adopts a probabilistic inventory control framework based on a Lévy process, which captures the stochastic nature of cumulative demand over time without relying on explicit forecasting. This formulation enables the derivation of closed-form expressions for expected total cost and facilitates analytical evaluation of how cost components are influenced by control parameters. In this way, the proposed approach offers complementary advantages to AI-based forecasting methods by emphasizing interpretability and mathematical tractability in decision-making under uncertainty.

\subsection{Prior Research on Inventory Management Using Stochastic Processes}

\noindent In the literature related to this study, the formulation of inventory models using stochastic processes, the optimality of $(s, S)$ policies, and analyses under discounted cost criteria have been widely investigated.  
Perera et al.~\cite{perera2023survey} provide a comprehensive survey of optimal strategies for stochastic inventory systems with fixed ordering costs in continuous-time settings.
\par In this context, the inventory level is often modeled as "initial inventory minus cumulative demand plus cumulative replenishment," as discussed in~\cite{muthuraman2015inventory},~\cite{bensoussan2006optimality}, and~\cite{bensoussan2010inventory}. These studies employ the probabilistic concept of stopping times to describe the timing of random replenishments.
\par In particular, Bensoussan et al.~\cite{bensoussan2010inventory} and Muthuraman et al.~\cite{muthuraman2015inventory} incorporate lead times into their formulations. However, in both studies, the timing of replenishments is determined independently of the cumulative demand size.
\par In such cases, when both demand and replenishment timing are random, it becomes difficult for inventory managers to issue replenishment orders appropriately, which may lead to risks of overstocking or stockouts.
\par Furthermore, Bensoussan et al.~\cite{bensoussan2006optimality} utilize conditional expectations and rely heavily on advanced probabilistic techniques. Although this approach preserves the information contained in the stochastic process within the mathematical model, it has the drawback that the formulation is based on discrete time and is difficult to extend to continuous time.
\par Recent advances in Lévy-process-based stochastic control have led to substantial progress in identifying optimal control strategies under complex conditions, including absolutely continuous control schemes~\cite{noba2025stochastic} and Poissonian intervention frameworks~\cite{yamazaki2017inventory}. While these studies provide elegant structural characterizations of optimality, they do not explicitly quantify the cost structure associated with key inventory control parameters such as reorder points and order quantities.
\par This study fills that gap by proposing an analytically tractable inventory model in which cumulative demand is modeled as a drifted Poisson process. The model enables closed-form approximations of the expected total cost and facilitates a detailed sensitivity analysis of how cost components are influenced by control parameters.

\section{Setting of Stochastic Models}

\noindent We consider the demand process $\{D_t\}_{t \ge 0}$ as a stochastic process, and define the cumulative demand as follows:

\begin{definition}[Cumulative Demand]
Let $\mu > 0$, $\alpha > 0$, and let $\{N_t\}_{t \ge 0}$ be a Poisson process with rate $\lambda > 0$. Then the cumulative demand process $\{D_t\}_{t \ge 0}$ is defined by
\begin{equation}
D_t = \mu t + \alpha N_t.
\end{equation}
\end{definition}

In this model, the inventory level at time $t$ is given by the initial inventory minus the cumulative demand up to time $t$. When the inventory level falls below a certain threshold, a replenishment order of fixed size $Q$ is placed. Since the times at which the inventory hits the threshold are random, we define corresponding stopping times.  
To model replenishment triggered at each reorder point, the system is formulated so that an order is placed whenever the cumulative demand exceeds a specified threshold, taking the order quantity into account.
\begin{definition}[First Passage Time (FPT)]
The first passage time $T_n$ for the process $\{D_t\}_{t \ge 0}$ is defined as
\begin{equation}
T_n = \inf \left\{s > 0 \,\middle|\, D_s \ge a + (n-1)Q \right\},
\end{equation}
where $a$ is the initial threshold  and $Q$ is the fixed order quantity.
\end{definition}
In general, it is known that the first passage time of a drifted Poisson process to a level $m$ approximately follows a gamma distribution $\Gamma(m, \lambda)$ under suitable scaling. Therefore, for the first passage time $T_n$ to the level $a + (n-1)Q$ in the process $D_t = \mu t + \alpha N_t$, we have the following approximation:

\begin{equation}
T_n \sim \Gamma\left( \frac{a + (n-1)Q}{\mu}, \alpha \lambda \right).
\end{equation}

The cumulative distribution function (CDF) of $T_n$ is then given by
\begin{equation}
P(T_n < t) = \int_0^t \frac{(\alpha \lambda s)^{\alpha(n)-1}}{(\alpha(n)-1)!} \, \alpha \lambda e^{-\alpha \lambda s} \, ds,
\end{equation}
where $t > 0$  and $\alpha(n) = \frac{a + (n-1)Q}{\mu}$.

\begin{definition}
The stochastic process $\Xt$ that represents the model where $x$ is the initial inventory ,$100-a$ is reorder point and $Q$ is ordered when the inventory quantity decreases to $x-a$ is defined as follows.
\begin{equation}
X_t = x - D_t + Q R_t .
\end{equation}
\end{definition}
Here,
\begin{equation}
R_t =  \sum_{n \ge 1}  1_{\{T_n <t \}} = \inf\{ k \mid T_k > t\}
\end{equation}
is called the number of renewals at $[0,t]$ with time $0$ as the renewal time.
\begin{figure}[H]
     \centering
     \includegraphics[width=0.8\textwidth]{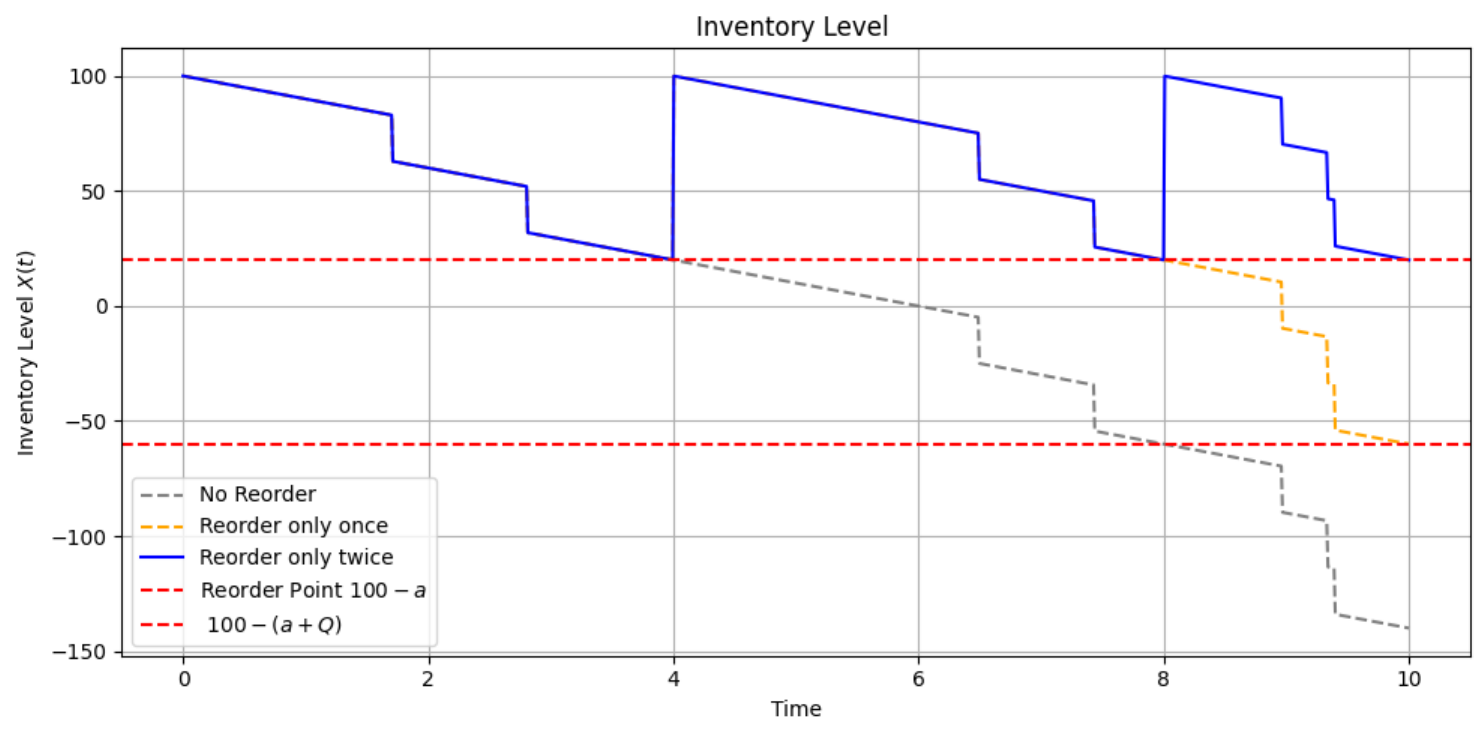}
     \caption{Inventory Level}
     \label{inventory level}
\end{figure}
In inventory control, “number of renewals” $=$ “number of orders placed,” so it can be used in simulations and in the calculation of ordering costs in total cost calculation.
 In fact, the expected quantity of inventory at time $t$ in this formulation can be calculated as follows.\begin{equation}
E[X_t] = x -\mu t -\alpha \lambda t +Q  \sum_{n \ge 1} \int_0^t \frac{(\alpha \lambda s)^{\alpha(n)-1}}{(\alpha(n)-1)!} s e^{-\alpha \lambda s}ds\\ .
\end{equation}

\begin{proof}
\begin{equation*}
    E[X_t] =  x - \mu t -\alpha \lambda t + Q E [R_t].
\end{equation*}
Also,by Fubini's Theorem,
\begin{align*}
E[R_t]
&=  \sum_{n \ge 1} E[1_{\{T_n< t\}}] =  \sum_{n \ge 1} P(T_n < t) =  \sum_{n \ge 1} \int_0^t \frac{(\alpha \lambda s)^{\alpha(n)-1}}{(\alpha(n)-1)!} s e^{-\alpha \lambda s}ds .\\
\end{align*}

\end{proof}

\section{Minimization of total cost}
\noindent In the previous section, inventory control was formulated as a stochastic process. Now we consider the problem of defining and minimizing the expected total cost.
\begin{definition}{(Expected total cost)}
Let $C_o$ denote the ordering cost per unit of inventory, $C_h$ the holding cost per unit per unit time, and $C_{s.o.}$ the stockout cost per unit per unit time.
The total cost per unit of inventory $TC(a,Q)$ is formulated as 
\begin{equation}
TC(a,Q,t) = \text{"ordering cost"} + \text{"holding cost"} + \text{"stockout cost"} .
\end{equation}
Define the expected value $E[TC(a,Q)]$ of the total cost $TC(a,Q)$ ,
\begin{equation}
E[TC(a,Q,t)] = Q C_o E [R_t] + C_h E\left[\int_0^t X_s ds\right] +C_{s.o.}E\left[\int_0^t (-X_s )^+ ds\right] .
\end{equation}
\end{definition}
From the way the model is built, we can assume that the likelihood of running out of stock is very small, and we can calculate
\begin{equation}
C_{s.o.}E\left[\int_0^t (-X_s )^+ ds\right] = \ C_{s.o.}E\left[\int_0^t (D_s - x -QR_s)^+ ds\right] \simeq 0 . \label{add}
\end{equation}
Applying the equation defined so far, we get
\begin{align*}
&\quad E[TC(a,Q,t)] \\ 
&= \ C_o Q E[R_t] + C_h E\left[\int_0^t(x - D_s + Q R_s) ds\right]   \\
&= \ C_o Q E[R_t] + C_h xt - C_h \int_0^t E [D_s]ds + C_h Q E \left[ \int_0^t R_sds \right]\\
&= \ C_o Q E[R_t] + C_h xt - C_h \int_0^t E [D_s]ds + C_h Q E \left[ \int_0^t R_sds \right]\\
&= \ C_o Q E[R_t] + C_h xt - C_h \int_0^t (\mu s+ \alpha \lambda s )ds + C_h Q E \left[ \int_0^t R_sds \right] \displaybreak[1]\\
&= \ C_o Q E[R_t] + C_h xt - C_h\frac{t^2}{2} (\mu + \alpha \lambda  ) + C_h Q E \left[ \int_0^t R_sds \right] . \\
\end{align*}
For the last term, we will use
\begin{align*}
E \left[\int_0^t R_s ds\right] 
&= \ E \left[\int_0^t \sum_{n \ge 1} 1_{\{T_n < s\}} ds\right] = \ \sum_{n \ge 1}E \left[\int_0^t  1_{\{T_n < s\}} ds\right] \\
&= \ \sum_{n \ge 1}E \left[(t -T_n) 1_{\{T_n < t\}}\right] = \ \sum_{n \ge 1} (E[t 1_{\{T_n < t\}}]-E[T_n1_{\{T_n < t\}}]) \\
&= \ \sum_{n \ge 1} (E[t 1_{\{T_n < t\}}]-E[T_n1_{\{T_n < t\}}]) \\
&= \ \sum_{n \ge 1} (t P(T_n < t)-E[T_n1_{\{T_n < t\}}]) \\
&= \ \sum_{n \ge 1} (t\int_0^t \frac{(\alpha \lambda s)^{\alpha(n)-1}}{(\alpha(n)-1)!} s e^{-\alpha \lambda s}ds
                       -\int_0^t s \frac{(\alpha \lambda s)^{\alpha(n)-1}}{(\alpha(n)-1)!} s e^{-\alpha \lambda s}ds ) .\\
\end{align*}
Based on the above calculations, the following approximate expression can be derived.
\begin{theorem}
Let $x > 0$ be the initial inventory level, $x - a > 0$ the reorder point, and $Q$ the order quantity. Let $C_o$, $C_h$, and $C_{s.o.}$ denote the ordering cost, holding cost, and shortage cost, respectively.  
Furthermore, suppose that the cumulative demand follows a drifted Poisson process with drift $\mu > 0$ and rate parameter $\lambda$, that is,
$D_t = \mu t + \alpha N_t$
,where $N_t$ is a Poisson process.
If the probability of stockout is extremely low, the expected total cost $E[TC(a, Q, t)]$ in this inventory control model can be approximately expressed as follows:
\begin{align} 
& \quad E[TC(a,Q,t)]  \notag \\
&\simeq \ C_o Q \sum_{n \ge 1} \int_0^t \frac{(\alpha \lambda s)^{\alpha(n)-1}}{(\alpha(n)-1)!} s e^{-\alpha \lambda s}ds + C_h xt - \frac{t^2}{2}C_h (\mu + \alpha \lambda  )  \notag \\
&+ C_h Q \sum_{n \ge 1} \left(t\int_0^t \frac{(\alpha \lambda s)^{\alpha(n)-1}}{(\alpha(n)-1)!} s e^{-\alpha \lambda s}ds 
-\int_0^t s \frac{(\alpha \lambda s)^{\alpha(n)-1}}{(\alpha(n)-1)!} s e^{-\alpha \lambda s}ds \right) \notag \\
\end{align}
\end{theorem}
,where $t>0$, $ \alpha(n) = \frac{a +(n-1) \times Q}{\mu} $ .
\section{Numerical results}

\noindent In this section, this studydescribe the simulation based on the expected total cost formulation introduced above. 
This studyset the initial inventory level to $x = 100$, and simulate the model using a drifted Poisson process with parameters $\mu = 5$, $\alpha = 10$, and $\lambda = 1$.  
Holding cost $C_h$, ordering cost $C_o$, shortage cost $C_{s.o.} $are selected such that they satisfy the relationship $C_h \le C_o \le C_{s.o.}$.
This studyinvestigate the effect of varying the reorder point a and order quantity Q on the total cost.
\par The risk of stockouts is extremely low, owing to the structure of our model. 
Therefore, the term involving shortage cost $C_{s.o.} $ can be effectively considered zero, which allows us to obtain an upper bound on the total cost. 
Generally, the expected value of a jump process smoothens individual jumps as time progresses, leading to a loss of information regarding jump sizes and their distribution. 
Consequently, it is not reliable to fully trust the global behavior of expectations when analyzing the underlying stochastic process. 
It is also noteworthy that the time at which the expected total cost reaches its maximum closely aligns with the convergence time of the numerical computation involving the gamma function. 
Thus, the expected number of replenishments can be accurately obtained around that point in time; however, beyond this, information about the jump sizes and their distribution is essentially lost. 
Consequently, the amount of useful information decreases, and it is reasonable to consider that these later periods are not particularly relevant for further analysis.
\par
\begin{figure}[H]
     \centering
     \includegraphics[width=0.8\textwidth]{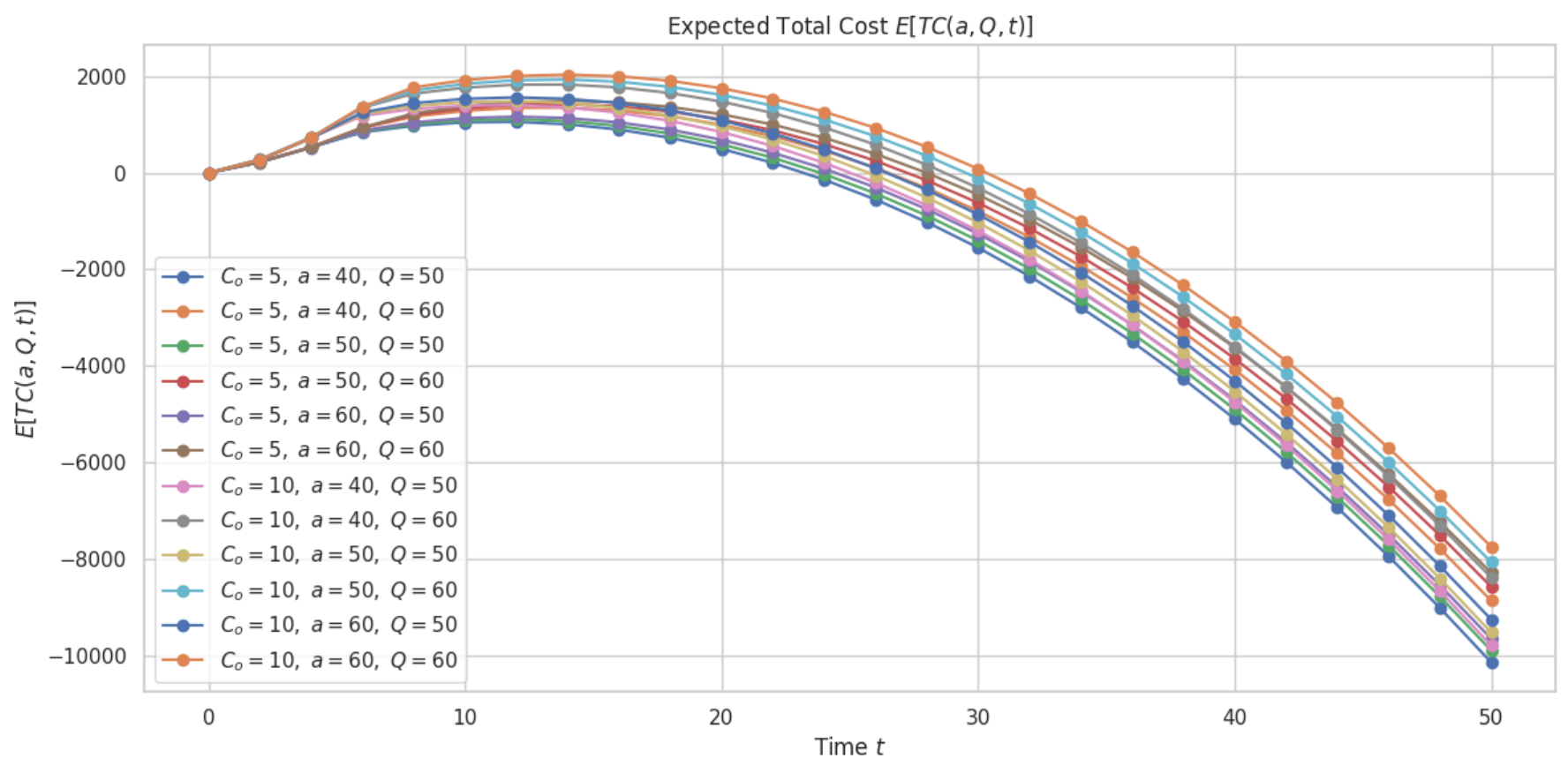}
     \caption{Expexted Total Cost}
     \label{totalcost}
\end{figure}

 The simulation results based on the formulation described in the previous section are shown in Figure~\ref{totalcost}.

\par The expected total cost of the proposed method was compared with the total cost obtained through a simulation using demand forecasts generated by an ARIMA model, followed by a reorder point policy. 
The comparison was performed according to the following procedure:
\par First, $1000 $ time series samples were generated. Because our proposed model involves a stochastic simulation, a direct comparison with the ARIMA model is not feasible. 
Therefore, this studygenerated 1000 sample paths from the drifted Poisson process, and performed simulations using the ARIMA model for each sample. 
By comparing the average results, this studyaimed to compare the expected value (i.e., the probabilistic mean) of the proposed model with the sample mean obtained from the ARIMA-based simulations.
\par Next, demand forecasting was performed for each time series using the ARIMA model. 
The ARIMA model has been a commonly used linear model for time series forecasting over the past $30$ years (Zhang~\cite{Zhang2003}). 
In addition, this studyemployed the rolling forecast method, which is widely used in demand–supply operations  (Huang~\cite{huang2011ordering}) , with a rolling window size of $12 $periods (i.e., one year). 
The parameters of the ARIMA model were tuned using the  \texttt{auto\_arima} function in the \texttt{pmdarima} Python library.  
 The search range for the ARIMA parameters was set to $p, q = 0, \ldots, 5$, and the differencing order d was automatically selected to ensure stationarity.
 \par Finally, an inventory simulation using the reorder point policy was performed for periods $13 - 50$, during which demand forecasting was applied. 
The initial inventory level, ordering cost, shortage cost, and holding cost were set to be the same as those used in the proposed method.
The total cost results from the simulations using the ARIMA model are presented in a Table~\ref{tab:avg-total-cost} at the end of this section.
\par
\begin{figure}[H]
     \centering
     \includegraphics[width=0.7\textwidth]{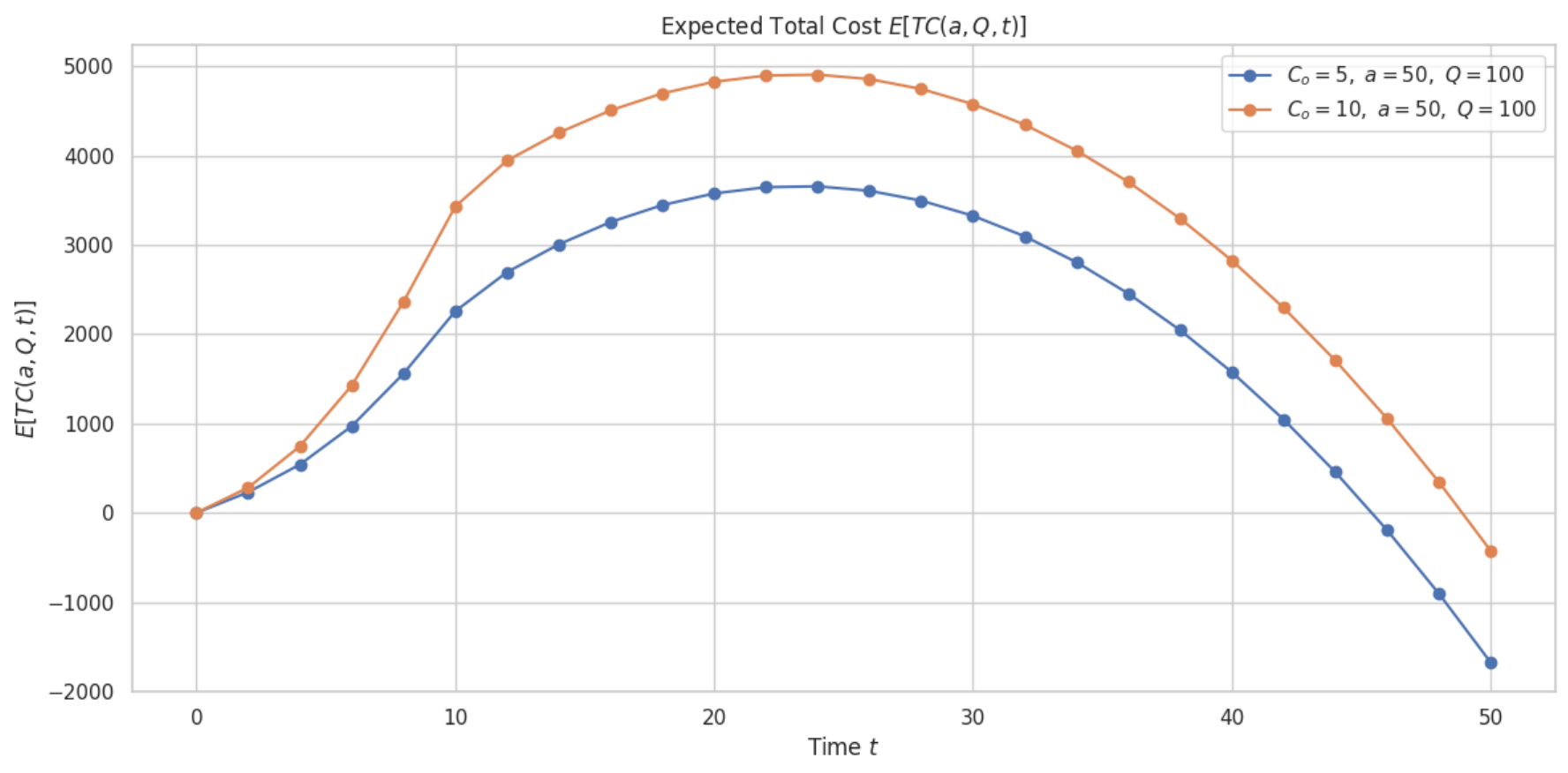}
     \caption{Fixed $a,Q$ }
     \label{fig:inventory1}
\end{figure}

 This studyanalyzed the expected total cost simulation by fixing specific parameters and individually examining their effects. 
First, this studyfixed the reorder point and order quantity, and varied the ordering cost. The results are shown in Figure~\ref{fig:inventory1}.
As the ordering cost increases, the total cost increases, which is consistent with our hypothesis.
\par
\begin{figure}[H]
     \centering
     \includegraphics[width=0.8\textwidth]{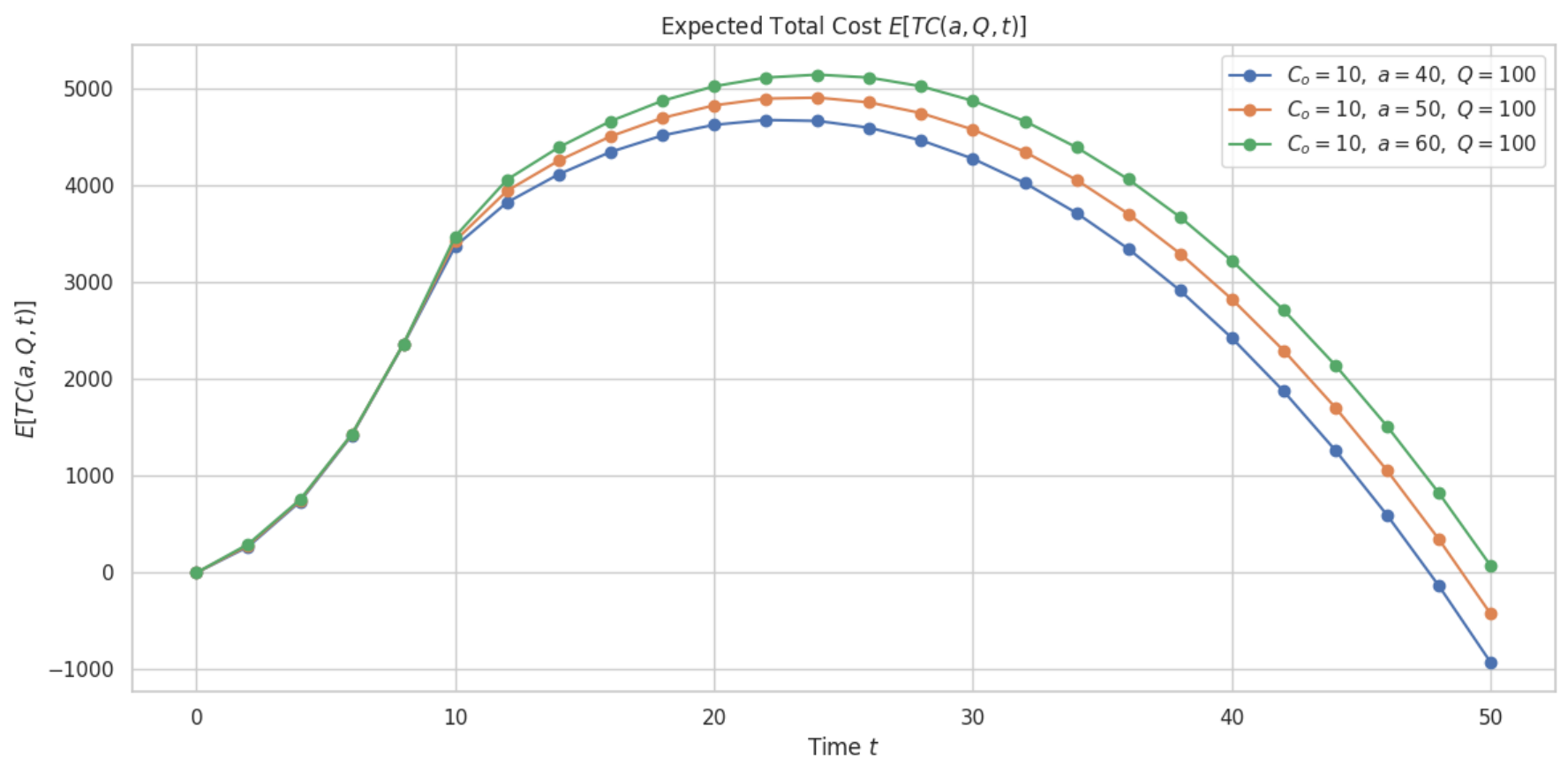}
     \caption{Fixed $Q,C_{o}$}
     \label{fig:inventory2}
\end{figure}

 Figure~\ref{fig:inventory2} shows the results when the ordering cost and order quantity are fixed, and the reorder point is varied. 
In this case, lowering the reorder point  (i.e. increasing $a$) decreases the number of orders, thereby reducing the total cost.
\par
\begin{figure}[H]
     \centering
     \includegraphics[width=0.8\textwidth]{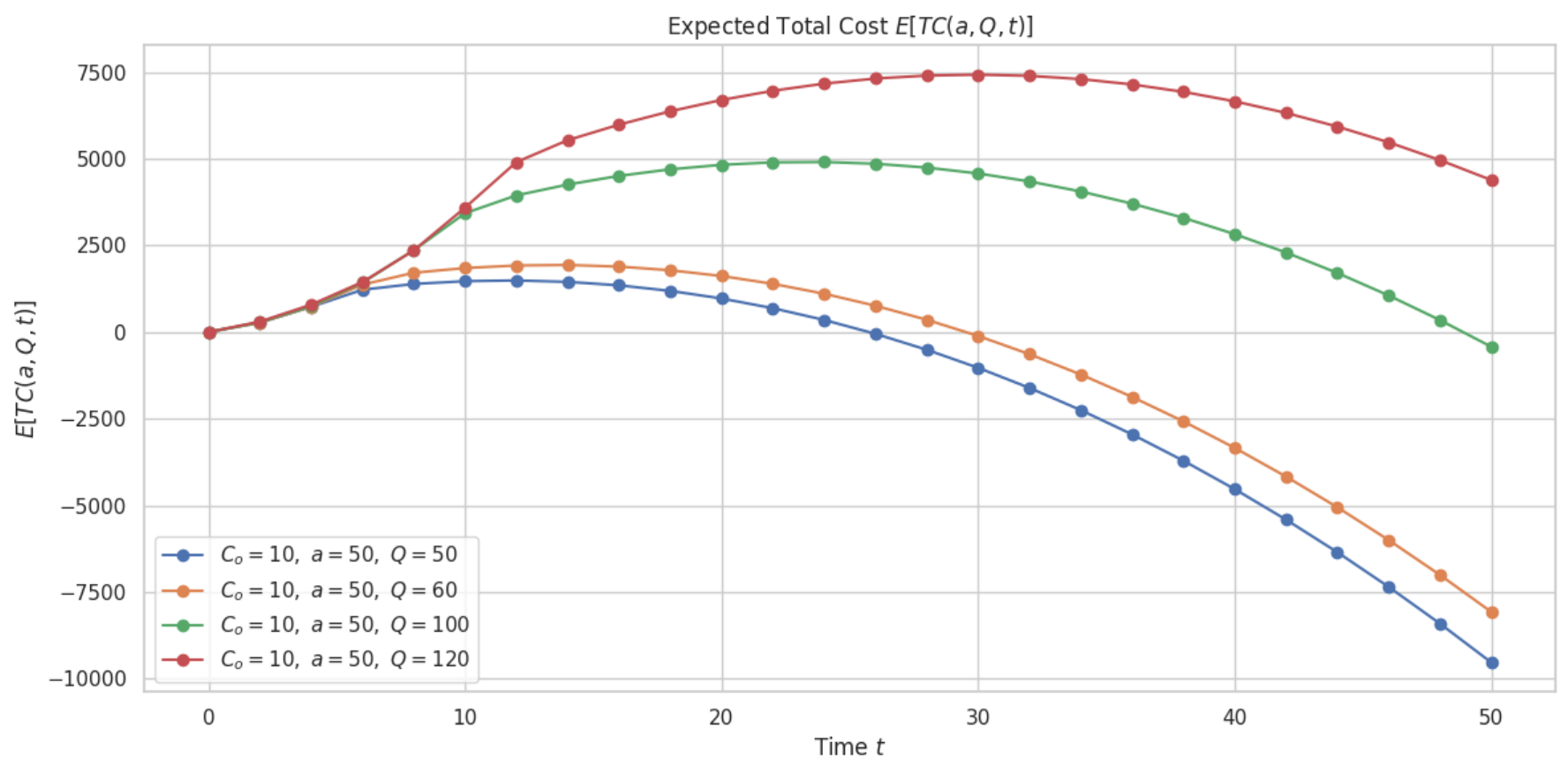}
     \caption{Fixed $a,C_{o}$}
     \label{fig:inventory3}
\end{figure}

 Finally, Figure~\ref{fig:inventory3} presents the results when the ordering cost and reorder point are fixed and the order quantity is varied. 
It can be observed that changes in the order quantity have a significant impact on the time-dependent behavior of the expected total cost.
\par It can be observed that the values of $t$ around the maximum in all these graphs are close to the point where the integrand of the gamma function in the proposed model assumes its dominant value. 
As the information related to jumps in a Lévy process is averaged over time, it is sufficient to analyze the graphs up to the vicinity of their maximum values.
\par For instance, this studyexamined the term-by-term values of the gamma function series used in the expectation calculation for  $a=50$ and $Q=60$. 
It was confirmed that the gamma function values converged sufficiently within a timeframe close to the cost function’s maximum.
\par
\begin{figure}[H]
     \centering
     \includegraphics[width=0.8\textwidth]{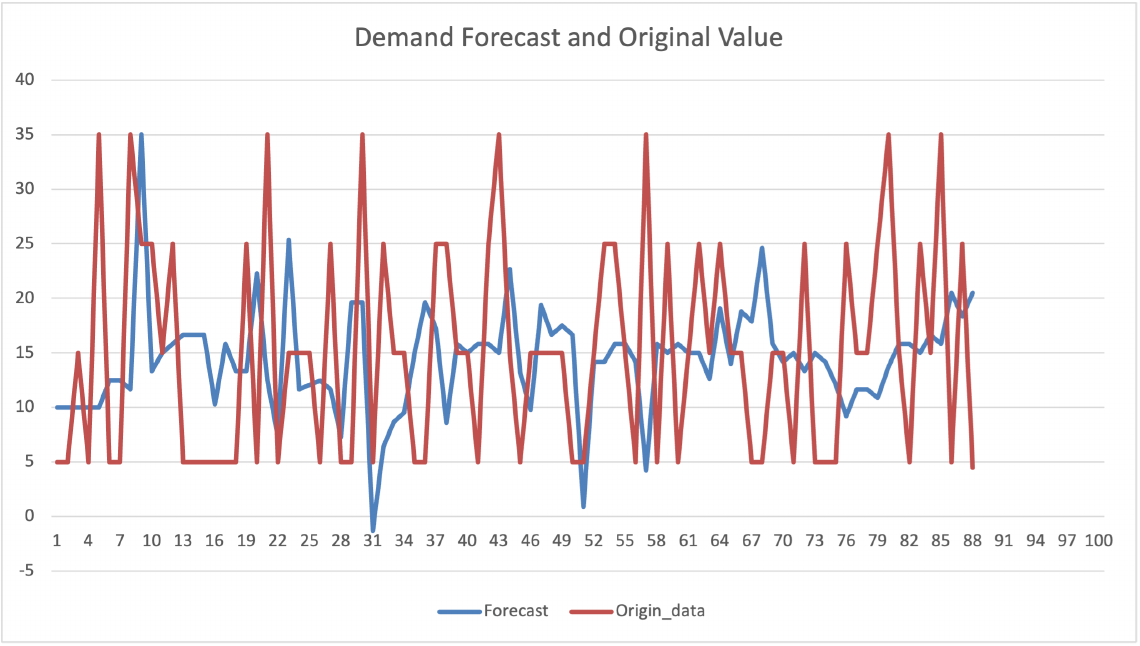}
     \caption{Demand Forecast and Original Value}
     \label{foracastarima}
\end{figure}

\begin{figure}
     \centering
     \includegraphics[width=0.8\textwidth]{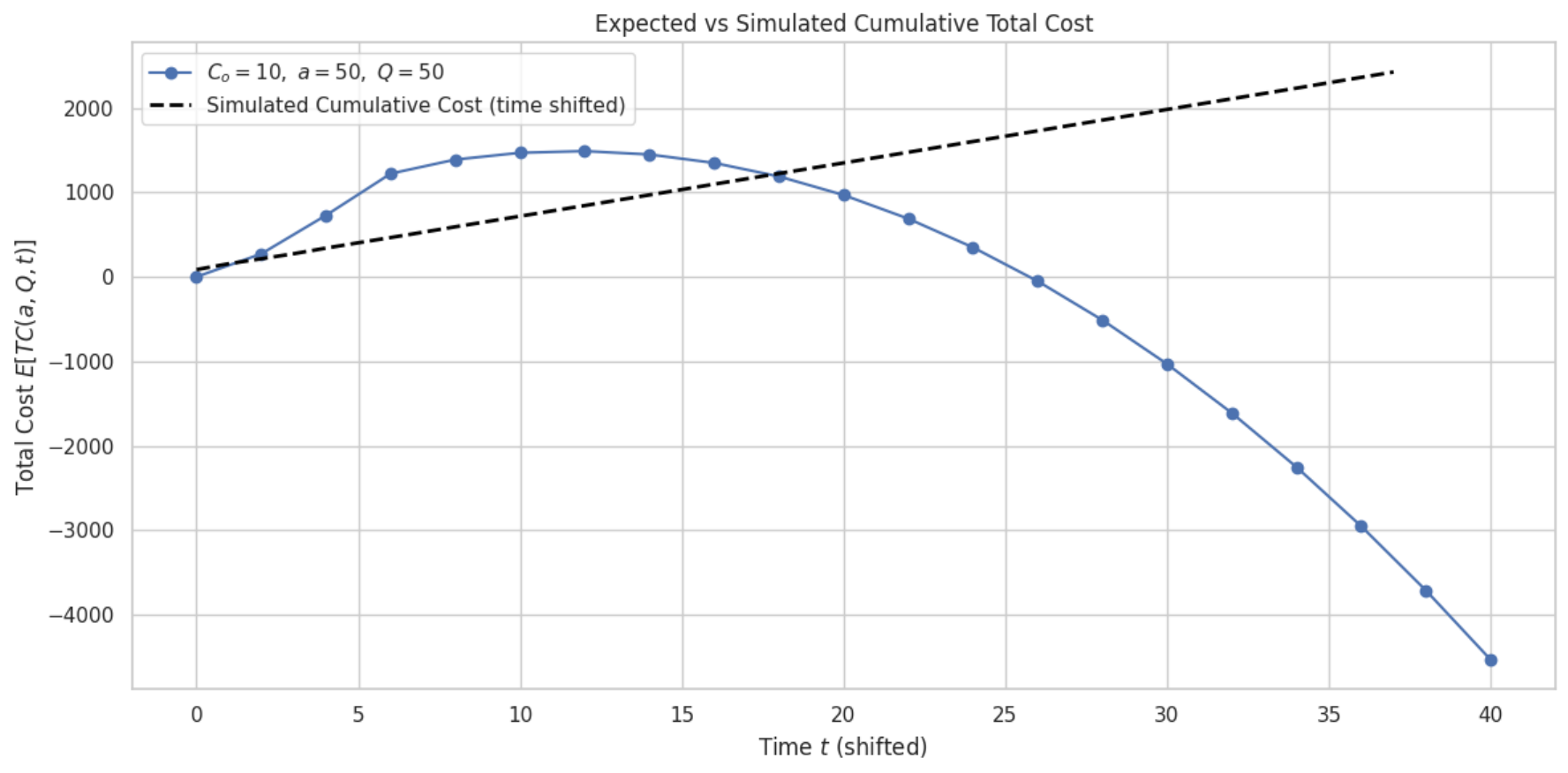}
     \caption{Expected Total cost vs Arima ($a=50,Q=50$)}
     \label{vsarima}
\end{figure}

\begin{figure}
     \centering
     \includegraphics[width=0.8\textwidth]{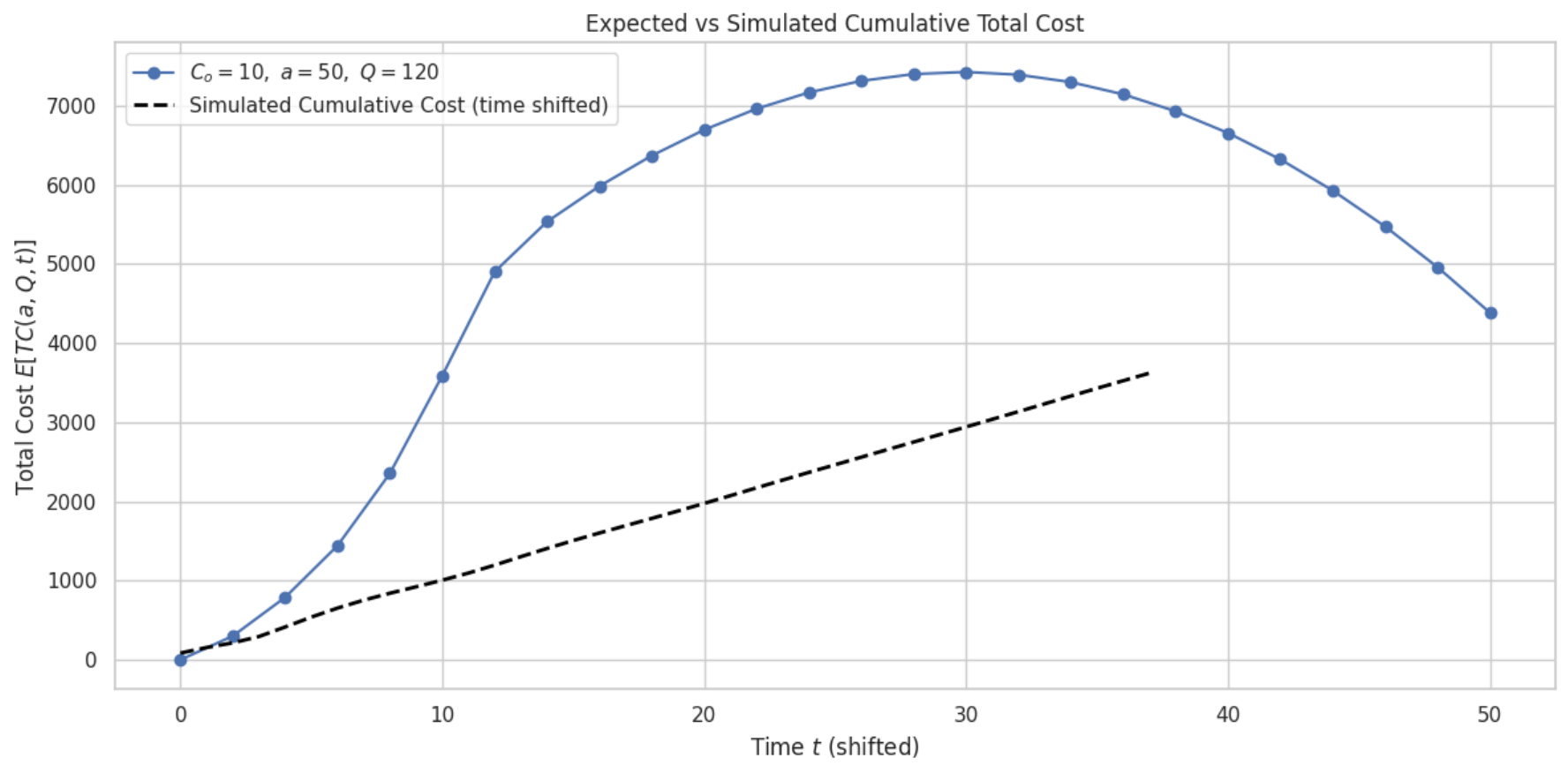}
     \caption{Expected Total cost vs Arima ($a=50,Q=120$)}
     \label{vsarimar50q120}
\end{figure}

 In addition, Figures~\ref{vsarima} and~\ref{vsarimar50q120} compare the total cost obtained from the proposed method with that of the ARIMA-based simulation using the same parameters. 
After training the ARIMA model for $12 $ periods, forecasting was performed; however, the model failed to accurately capture the fluctuations in the original time series. 
Specifically, the ARIMA model tended to underestimate future demand (Figures~\ref{foracastarima}).
\begin{figure}[H]
     \centering
     \includegraphics[width=0.8\textwidth]{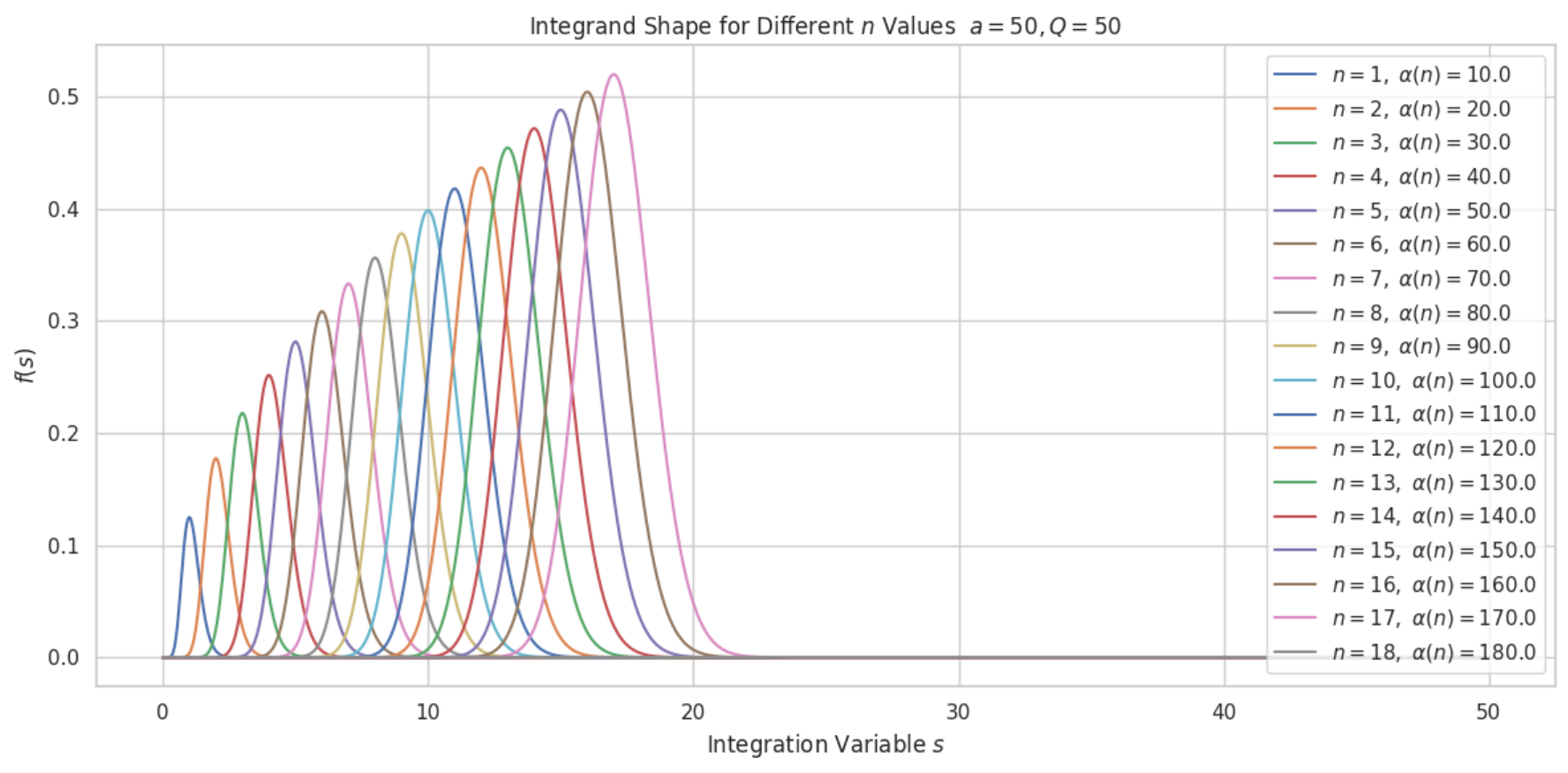}
     \caption{Gamma Function}
     \label{gamma}
\end{figure}

\par Considering this, it is argued that our expected total cost, which exceeds the cost obtained from the ARIMA-based simulation, should not be viewed as an overestimation, but rather as a reasonable and justified result.

\begin{table}[H]
  \centering
  \caption{Average Total Cost by Parameter Set}
  \label{tab:avg-total-cost}
  \small
  \setlength{\tabcolsep}{4pt}
  \begin{minipage}[t]{0.48\textwidth}
    \centering
    \begin{tabular}{rrrrrr}
      \toprule
      $R$ & $Q$ & $C_h$ & $C_o$ & $C_{s.o.}$ & Total Cost \\
      \midrule
      40 & 50 & 1 & 5 & 10 & 2108 \\
      40 & 50 & 1 & 10 & 10 & 2161 \\
      40 & 50 & 1 & 5 & 15 & 2157 \\
      40 & 50 & 1 & 10 & 15 & 2210 \\
      40 & 60 & 1 & 5 & 10 & 2269 \\
      40 & 60 & 1 & 10 & 10 & 2314 \\
      40 & 60 & 1 & 5 & 15 & 2316 \\
      40 & 60 & 1 & 10 & 15 & 2361 \\
      50 & 50 & 1 & 5 & 10 & 2464 \\
      50 & 50 & 1 & 10 & 10 & 2518 \\
      50 & 50 & 1 & 5 & 15 & 2512 \\
      50 & 50 & 1 & 10 & 15 & 2567 \\
      50 & 60 & 1 & 5 & 10 & 2626 \\
      50 & 60 & 1 & 10 & 10 & 2672 \\
      50 & 60 & 1 & 5 & 15 & 2672 \\
      50 & 60 & 1 & 10 & 15 & 2718 \\
      60 & 50 & 1 & 5 & 10 & 2823 \\
      60 & 50 & 1 & 10 & 10 & 2879 \\
      60 & 50 & 1 & 5 & 15 & 2872 \\
      60 & 50 & 1 & 10 & 15 & 2927 \\
      60 & 60 & 1 & 5 & 10 & 2990 \\
      60 & 60 & 1 & 10 & 10 & 3037 \\
      60 & 60 & 1 & 5 & 15 & 3037 \\
      60 & 60 & 1 & 10 & 15 & 3083 \\
      \bottomrule
    \end{tabular}
  \end{minipage}
  \begin{minipage}[t]{0.48\textwidth}
    \centering
    \begin{tabular}{rrrrrr}
      \toprule
      $R$ & $Q$ & $C_h$ & $C_o$ & $C_{s.o.}$ & Total Cost \\
      \midrule
      40 & 110 & 1 & 5 & 10 & 3110 \\
      40 & 110 & 1 & 10 & 10 & 3136 \\
      40 & 110 & 1 & 5 & 15 & 3147 \\
      40 & 110 & 1 & 10 & 15 & 3173 \\
      40 & 120 & 1 & 5 & 10 & 3288 \\
      40 & 120 & 1 & 10 & 10 & 3312 \\
      40 & 120 & 1 & 5 & 15 & 3324 \\
      40 & 120 & 1 & 10 & 15 & 3348 \\
      50 & 110 & 1 & 5 & 10 & 3485 \\
      50 & 110 & 1 & 10 & 10 & 3511 \\
      50 & 110 & 1 & 5 & 15 & 3522 \\
      50 & 110 & 1 & 10 & 15 & 3548 \\
      50 & 120 & 1 & 5 & 10 & 3665 \\
      50 & 120 & 1 & 10 & 10 & 3690 \\
      50 & 120 & 1 & 5 & 15 & 3701 \\
      50 & 120 & 1 & 10 & 15 & 3725 \\
      60 & 110 & 1 & 5 & 10 & 3865 \\
      60 & 110 & 1 & 10 & 10 & 3892 \\
      60 & 110 & 1 & 5 & 15 & 3903 \\
      60 & 110 & 1 & 10 & 15 & 3929 \\
      60 & 120 & 1 & 5 & 10 & 4048 \\
      60 & 120 & 1 & 10 & 10 & 4073 \\
      60 & 120 & 1 & 5 & 15 & 4084 \\
      60 & 120 & 1 & 10 & 15 & 4108 \\
      \bottomrule
    \end{tabular}
  \end{minipage}
\end{table}

\section{Conclusion}
\noindent
This study has proposed an analytical model for inventory control under intermittent demand, formulating cumulative demand as a drifted Poisson process. 
By explicitly modeling both demand and replenishment as stochastic processes, we derived a closed-form expression for the expected total cost using tools from stochastic process theory.
\par
The resulting cost function enables sensitivity analysis with respect to key parameters, such as order quantity, reorder point, and demand characteristics, providing deeper insight into their influence on the overall cost structure.
Owing to its analytical tractability, the model offers notable computational advantages over simulation-based and machine learning approaches, such as neural network forecasting methods~\cite{KOURENTZES2013198}, making it well-suited for real-time decision support.
\par
Future research will extend this framework to incorporate statistical risk management metrics, including stockout probabilities and cost variability, thereby bridging practical decision-making with theoretical rigor. In particular, we aim to explore inventory systems driven by Lévy processes, which can better capture demand shocks and heavy-tailed behaviors~\cite{yamazaki2017inventory}. 
This extension will leverage the Lévy–Khintchine representation~\cite{sato1999levy}~\cite{applebaum2009levy} and stochastic delay differential equations (SDDEs)~\cite{OksendalSulem2007} to more flexibly model information lags and stochastic lead times.
\par
In particular, our previous study~\cite{Koide2024} formulated a framework for analyzing the impact of time delays in inventory adjustment using Lévy processes. While such formulations offer expressive modeling power, applying SDDEs requires careful consideration of the theoretical properties of the solution, especially the conditions for existence and uniqueness.
\par
Moreover, integrating stochastic control techniques and established optimality results under information delays~\cite{bensoussan2006optimality}~\cite{bensoussan2010inventory} will be critical for constructing a practical and theoretically sound framework for inventory control under uncertainty.

\bibliographystyle{plain}
\bibliography{ref}

\begin{thebibliography}{10}

\bibitem{altay2008adapting}
Nezih Altay, Frank Rudisill, and Lewis~A Litteral.
\newblock Adapting wright's modification of holt's method to forecasting
  intermittent demand.
\newblock {\em International Journal of Production Economics}, 111(2):389--408,
  2008.

\bibitem{applebaum2009levy}
David Applebaum.
\newblock {\em L{\'e}vy Processes and Stochastic Calculus}, volume 116.
\newblock Cambridge university press, 2009.

\bibitem{bensoussan2006optimality}
Alain Bensoussan, Metin {\c{C}}akanyildirim, and Suresh~P Sethi.
\newblock Optimality of base-stock and (s, s) policies for inventory problems
  with information delays.
\newblock {\em Journal of optimization theory and applications}, 130:153--172,
  2006.

\bibitem{bensoussan2010inventory}
Alain Bensoussan, Lama Moussawi-Haidar, and Metin {\c{C}}akany{\i}ld{\i}r{\i}m.
\newblock Inventory control with an order-time constraint: optimality,
  uniqueness and significance.
\newblock {\em Annals of Operations Research}, 181:603--640, 2010.

\bibitem{croston1972forecasting}
J.D. Croston.
\newblock Forecasting and stock control for intermittent demands.
\newblock {\em Operational Research Quarterly}, 23(3):289--303, 1972.

\bibitem{huang2011ordering}
Li-Ting Huang, I-Chien Hsieh, and Cheng-Kiang Farn.
\newblock On ordering adjustment policy under rolling forecast in supply chain
  planning.
\newblock {\em Computers \& Industrial Engineering}, 60(3):397--410, 2011.

\bibitem{Koide2024}
Ryoya Koide, Yurika Ono, and Aya Ishigaki.
\newblock Impact analysis on the time delay of inventory adjustment using the
  {L}\'evy process.
\newblock {\em Proceedings of the Japan Joint Automatic Control Conference},
  67:235--236, 2024.

\bibitem{KOURENTZES2013198}
Nikolaos Kourentzes.
\newblock Intermittent demand forecasts with neural networks.
\newblock {\em International Journal of Production Economics}, 143(1):198--206,
  2013.

\bibitem{muthuraman2015inventory}
Kumar Muthuraman, Sridhar Seshadri, and Qi~Wu.
\newblock Inventory management with stochastic lead times.
\newblock {\em Mathematics of Operations Research}, 40(2):302--327, 2015.

\bibitem{MR4702096}
Kei Noba and Kazutoshi Yamazaki.
\newblock On singular control for {L}\'evy processes.
\newblock {\em Mathematics of Operations Research}, 48(3):1213--1234, 2023.

\bibitem{noba2025stochastic}
Kei Noba and Kazutoshi Yamazaki.
\newblock On stochastic control under poissonian intervention: optimality of a
  barrier strategy in a general l{\'e}vy model.
\newblock {\em Journal of Applied Probability}, pages 1--24, 2025.

\bibitem{OksendalSulem2007}
Bernt {\O}ksendal and Agn\`es Sulem.
\newblock {\em Applied Stochastic Control of Jump Diffusions}.
\newblock Springer, 2nd edition, 2007.

\bibitem{perera2023survey}
Sandun~C Perera and Suresh~P Sethi.
\newblock A survey of stochastic inventory models with fixed costs: Optimality
  of (s, s) and (s, s)-type policies—discrete-time case.
\newblock {\em Production and Operations Management}, 32(1):131--153, 2023.

\bibitem{sato1999levy}
{Ken-Iti} Sato.
\newblock {\em L{\'e}vy Processes and Infinitely Divisible Distributions},
  volume~68.
\newblock Cambridge university press, 1999.

\bibitem{WANG20241038}
Shengjie Wang, Yanfei Kang, and Fotios Petropoulos.
\newblock Combining probabilistic forecasts of intermittent demand.
\newblock {\em European Journal of Operational Research}, 315(3):1038--1048,
  2024.

\bibitem{yamazaki2017inventory}
Kazutoshi Yamazaki.
\newblock Inventory control for spectrally positive {L}{\'e}vy demand
  processes.
\newblock {\em Mathematics of Operations Research}, 42(1):212--237, 2017.

\bibitem{yuna2023inventory}
Ferhat Yuna, Burak Erkayman, and Mustafa Y{\i}lmaz.
\newblock Inventory control model for intermittent demand: a comparison of
  metaheuristics.
\newblock {\em Soft Computing}, 27(10):6487--6505, 2023.

\bibitem{Zhang2003}
G.~Peter Zhang.
\newblock Time series forecasting using a hybrid arima and neural network
  model.
\newblock {\em Neurocomputing}, 50:159--175, 2003.

\bibitem{Zhou2021}
Haoyi Zhou, Shanghang Zhang, Jieqi Peng, Shuai Zhang, Jianxin Li, Hui Xiong,
  and Wancai Zhang.
\newblock Informer: Beyond efficient transformer for long sequence time-series
  forecasting.
\newblock In {\em Proceedings of the AAAI Conference on Artificial
  Intelligence}, volume~35, pages 11106--11115, 2021.

\end{thebibliography}

\end{document}